\documentclass[12pt]{amsart}

\usepackage{amssymb}
\usepackage{amsthm}
\usepackage{doc}
\usepackage{latexsym}
\usepackage{amscd}
\usepackage{graphics}
\usepackage{epsfig}
\usepackage{eucal}
\usepackage{mathrsfs}
\numberwithin{equation}{section}
\title{A note on representations of Lie-Yamaguti algebras induced by left Leibniz algebras}

\author{A. Nourou ISSA}
\address{D\'epartement de Math\'ematiques, Universit\'e d'Abomey-Calavi,\\
01 BP 4521  Cotonou 01, Benin}
\email{woraniss@yahoo.fr}
\thanks{}

\date{}

\begin{document}

\begin{abstract}
 It is well-known that each left Leibniz algebra has a natural structure of a Lie-Yamaguti algebra. In this paper it is shown that every left representation of a left Leibniz algebra
 $(\mathfrak{g}, \cdot)$ induces naturally a representation of the Lie-Yamaguti algebra
 $(\mathfrak{g}, [,], [\![, ,]\!])$  that is associated with $(\mathfrak{g}, \cdot)$. Moreover, it is proved that equivalent representations of $(\mathfrak{g}, \cdot)$ give equivalent representations of
 $(\mathfrak{g}, [,], [\![, ,]\!])$.


\end{abstract}


\footnotetext{
{\small\it 2020 Subject Classification:} 17A32, 17B60, 17D99

{\small\it Key words and phrases:} Leibniz algebra, Lie-Yamaguti algebra, representation
}



\maketitle

\section{Introduction}

The notion of a ``$D$-algebra'' was introduced by A. M. Bloh (\cite{Bl}). About three decades later, developing a homology theory for Lie algebras, J. L. Loday (\cite{L1, L2}) rediscovered this kind of algebras and called it ``Leibniz algebra''. As a generalization of Lie algebras and because of their applications in mathematics, physics, and mathematical physics, Leibniz algebras were studied in last decades from various aspects by many researchers. A representation and cohomology theory for Leibniz algebras was initiated in \cite{LP}.
\par
The notion of a ``generalized Lie triple system'' was introduced by K. Yamaguti (\cite{Yam1}) who developed further a representation and cohomology theory for this kind of algebras (\cite{Yam2}, \cite{Yam3}). In \cite{KW} M. K. Kinyon and A. Weinstein renamed generalized Lie triple systems as ``Lie-Yamaguti algebras''; moreover, they found that any left Leibniz algebra has a natural Lie-Yamaguti structure. Later on, investigating the Akivis algebra associated with a left Leibniz algebra, this relationship with Lie-Yamaguti algebras was considered again in \cite{Iss1}. Thus some specific constructions regarding Lie-Yamaguti algebras can be deduced from the ones of Leibniz algebras.
\par
In this short note we consider the representation of a Lie-Yamaguti algebra that is induced by a left representation of a given left Leibniz algebra. Recall that any representation of a given Maltsev algebra $M$ induces a representation of the Lie-Yamaguti algebra that is associated with $M$ (\cite{Yam2}). In section 2 useful notions and results are reminded and the main result is proved. In section 3, starting from some well-known representations of left Leibniz algebras, our main result is applied to give a few examples of nontrivial representations of Lie-Yamaguti algebras. All vectors and algebras are considered over a field of characteristic zero.

\section{Main results}

In this section we first remind some useful basic notions and next we prove that any left representation of a given left Leibniz algebra
$(\mathfrak{g}, \cdot)$ induces a representation of the associated Lie-Yamaguti algebra; moreover, it turns out that equivalent left representations of $(\mathfrak{g}, \cdot)$ imply equivalent representations the associated Lie-Yamaguti algebra. \\

\noindent
{\bf Definition 2.1.} (\cite{Bl, L1, L2}) A (left) Leibniz algebra is a vector space
$\mathfrak{g}$ together with a bilinear map
$\cdot : \mathfrak{g} \times \mathfrak{g} \rightarrow \mathfrak{g}$ satisfying the {\it left Leibniz identity}
\par
 $x \cdot (y \cdot z) = (x \cdot y) \cdot z + y \cdot (x \cdot z)$\\
for all $x,y,z \in \mathfrak{g}$.\\
\par
Observe that in \cite{L1, L2} only {\it right} Leibniz algebras were considered, i.e. algebras satisfying the {\it right Leibniz identity}
$(z \cdot y) \cdot x = z \cdot (y \cdot x) + (z \cdot x) \cdot y$.\\

\noindent
{\bf Definition 2.2.} (\cite{L2, LP}) Let
$(\mathfrak{g}, \cdot)$ be a left Leibniz algebra. A left representation of
$(\mathfrak{g}, \cdot)$ is a triple $(V,l,r)$, where $V$ is a vector space, $l,r:\mathfrak{g} \rightarrow End(V)$ are linear maps such that the following conditions hold for all $x,y \in
 \mathfrak{g}$:
 \begin{equation}\label{eq2.1}
 l(x \cdot y) = [l(x),l(y)],
 \end{equation}
 \begin{equation}\label{eq2.2}
 r(x \cdot y) = r(y) \circ r(x) + l(x) \circ r(y),
 \end{equation}
 \begin{equation}\label{eq2.3}
 r(x \cdot y) = [l(x),r(y)].
 \end{equation}
The vector space $V$ is also called a $\mathfrak{g}$-module.\\
\par
The conditions (\ref{eq2.2}) and (\ref{eq2.3}) imply
\begin{equation}\label{eq2.4}
 r(y) \circ r(x) + r(y) \circ l(x) = 0.
\end{equation}
\par
In the sequel, given some maps $f$ and $g$, their usual composition $f \circ g$ will be denoted simply by $fg$.
\par
An example of a left representation of
$(\mathfrak{g}, \cdot)$ is the so-called {\it adjoint} (or {\it regular}) representation that is defined as the triple $(\mathfrak{g}, L, R)$ with linear maps $L$, $R$ defined respectively by $L : \mathfrak{g} \rightarrow End(V)$, $x \mapsto L_x$ and $R : \mathfrak{g} \rightarrow End(V)$, $x \mapsto R_x$, where $L_x$ (resp. $R_x$) denotes the left (resp. right) multiplication in ($\mathfrak{g}, \cdot$), $L_x y := x \cdot y$ and $R_x y := y \cdot x$.
\par
Let $V^*$ denotes the dual space of $V$ and define two linear maps $l^* : \mathfrak{g} \rightarrow End (V^*)$, $x \mapsto l^* (x)$ and $r^* : \mathfrak{g} \rightarrow End (V^*)$, $x \mapsto r^* (x)$, where $l^* (x) \xi := - \xi  l(x)$, and $r^* (x) \xi := - \xi r(x)$ for all $x \in \mathfrak{g}$, $\xi \in V^*$. It is well-known that, in general, ($V^* , l^*, r^*$) is not a left representation of ($\mathfrak{g}, \cdot$). However, the following result provides a left representation of ($\mathfrak{g}, \cdot$).\\

\noindent
{\bf Lemma 2.3.} (\cite{TS})
 {\it Let $(\mathfrak{g}, \cdot)$ be a left Leibniz algebra and $(V,l,r)$ its left representation. Then $(V^*, l^*, -l^*-r^*)$ is a left representation of
 $(\mathfrak{g}, \cdot)$}.
\par
The representation $(V^*, l^*, -l^*-r^*)$ is called the {\it dual representation} of
$(\mathfrak{g}, \cdot)$.\\

\noindent
{\bf Definition 2.4.} (\cite{Yam1})
 A Lie-Yamaguti algebra (LY algebra for short) is a triple $(\mathfrak{g}, [,], [\![, ,]\!])$, where
 $\mathfrak{g}$ is a vector space, $[,]: \mathfrak{g} \times \mathfrak{g} \rightarrow \mathfrak{g}$ a bilinear map, and $[\![ , , ]\!]: \mathfrak{g} \times \mathfrak{g} \times \mathfrak{g} \rightarrow \mathfrak{g}$ a trilinear map such that
 \par
 (LY01) $[x,y] = - [y,x]$,
 \par
 (LY02) $[\![ x,y,z ]\!] = -[\![ y,x,z ]\!]$,
 \par
 (LY1) ${\circlearrowleft}_{x,y,z} ([[x,y],z] + [\![ x,y,z ]\!]) = 0$,
 \par
 (LY2) ${\circlearrowleft}_{x,y,z} ([\![[x,y],z,u ]\!]) = 0$,
 \par
 (LY3) $[\![ x,y,[u,v] ]\!] = [ [\![ x,y,u ]\!],v] + [u, [\![ x,y,v ]\!]]$,
 \par
 (LY4) $[\![ x,y, [\![ u,v,w ]\!] ]\!]=[\![ [\![x,y,u]\!], v,w ]\!] + [\![ u,[\![x,y,v]\!],w ]\!] + [\![ u,v,[\![x,y,w]\!] ]\!]$\\
 for all $u,v,w,x,y,z \in \mathfrak{g}$, where ${\circlearrowleft}_{x,y,z}$ denotes the sum over cyclic permutation of $x,y,z$.\\
\par
In \cite{KW} a relationship between left Leibniz algebras and $LY$ algebras has been discovered. Specifically, any left Leibniz algebra
$(\mathfrak{g}, \cdot)$ has a natural $LY$ structure (the associate $LY$ algebra) with respect to the operations
\par
$[x,y] := x \cdot y - y \cdot x$,
\par
$[\![x,y,z]\!] := - (x \cdot y) \cdot z$\\
for all $x,y,z \in \mathfrak{g}$ (see also \cite{Iss1}). Thus, examples of $LY$ algebras (even with specific properties) could be constructed from Leibniz algebras (e.g. see \cite{Iss2}).
\par
A representation and cohomology theory for $LY$ algebras was initiated in \cite{Yam2}, \cite{Yam3}. Their $(2,3)$-cohomology was applied in \cite{ZL} to study their deformations and abelian extensions.\\

\noindent
{\bf Definition 2.5.} (\cite{Yam2}) Let
$(\mathfrak{g}, [,], [\![, ,]\!])$ be a $LY$ algebra. A representation of $(\mathfrak{g}, [,], [\![, ,]\!])$ is a quadruple $(V, \rho, \theta, D)$ where $V$ is a vector space, $\rho : \mathfrak{g} \rightarrow End(V)$ is a linear map, and $\theta , D : \mathfrak{g} \times \mathfrak{g} \rightarrow End(V)$ are bilinear maps such that
\par
(R1) $D(x,y) - \theta (y, x) + \theta (x, y) + \rho ([x, y]) - [\rho (x), \rho (y)] = 0$;
\par
(R2) $D([x, y], z) + D([y, z], x) + D([z, x], y) =0$;
\par
(R3) $\theta ([x, y], z) = \theta (x, z) \rho (y) - \theta (y, z) \rho (x)$;
\par
(R4) $[D(x, y), \rho (z)] = \rho ([\![x,y,z]\!])$;
\par
(R5) $\theta (x, [y, z]) = \rho (y) \theta (x, z)-\rho (z) \theta (x, y)$;
\par
(R6) $[D(u, v), \theta (x, y)] = \theta ([\![u, v, x]\!], y) + \theta (x,
[\![u, v, y]\!])$;
\par
(R7) $\theta (u, [\![x, y, z]\!]) =
\theta (y, z)\theta (u, x) - \theta (x, z)\theta (u, y) + D(x, y)\theta (u, z)$\\
for all $u, v, x, y, z \in \mathfrak{g}$.\\
\par
An example of a representation of
$(\mathfrak{g}, [,], [\![, ,]\!])$ is the so-called {\it adjoint} (or {\it regular}) representation  $(\mathfrak{g}, \rho, \theta, D)$, where $\rho (x)(y) := [x,y]$, $\theta (x,y)(z) := [\![z,x ,y]\!]$, and $D (x,y)(z):= [\![x,y ,z]\!]$ for all $x,y,z \in \mathfrak{g}$ (\cite{Yam3}).
\par
Because of (R1), the representation $(V, \rho, \theta, D)$ is simply written as $(V, \rho, \theta)$. Our main result is as follows.\\

\noindent
{\bf Theorem 2.6.} 
{\it Let  $(\mathfrak{g}, \cdot)$ be a Leibniz algebra and $(V,l,r)$ its left representation. Let $(\mathfrak{g}, [,], [\![, ,]\!])$ be the $LY$ algebra associated with $(\mathfrak{g}, \cdot)$. Then $(V,l,r)$ induces a representation of $(\mathfrak{g}, [,], [\![, ,]\!])$ (call it the associated representation).}
\begin{proof}
 On $(\mathfrak{g}, \cdot)$ we define the following maps, for all $x,y \in \mathfrak{g}$:
\par
$\rho (x) := l(x) - r(x)$;
\par
$\theta (x,y) := -r(y)r(x)$;
\par
$D(x,y) := -l(x \cdot y)$. \\
Clearly $\rho : \mathfrak{g} \rightarrow End(V)$ is a linear map and $\theta , D : \mathfrak{g} \times \mathfrak{g} \rightarrow End(V)$ are bilinear maps. We proceed to check that the triple $(\rho, \theta, D)$ as defined above satisfies the set of axioms as in Definition 2.5.
\par
$\bullet$ For (R1) we have
\par
$D(x,y) + \theta (x, y) - \theta (y,x) + \rho ([x, y]) - [\rho (x), \rho (y)]$
\par
$= \{ -l(y \cdot x) + l(y)l(x) - l(x)l(y) \}_{(2.1)}$
\par
$+ \{- r(x \cdot y) + l(x)r(y) - r(y)l(x) \}_{(\ref{eq2.3})}$
\par
$+ \{r(y \cdot x) - l(y)r(x) + r(x)l(y) \}_{(\ref{eq2.3})}$
\par
$= 0$ (by (\ref{eq2.1}) and (\ref{eq2.3}))\\
so (R1) holds for $(\rho, \theta, D)$.
\par
$\bullet$ For (R2) we have
\par
$D([x,y],z) + D([y,z],x) + D([z,x],y)$
\par
$= - l([x,y] \cdot z) - l([y,z] \cdot x) - l([z,x] \cdot y)$
\par
$= - (l([x,y])l(z) - l(z)l([x,y])) - (l([y,z])l(x) - l(x)l([y,z])) - (l([z,x])l(y) - l(y)l([z,x]))$ (by (\ref{eq2.1}))
\par
$= - l(x \cdot y)l(z) + l(y \cdot x)l(z) + l(z)l(x \cdot y) - l(z)l(y \cdot x) - l(y \cdot z)l(x) + l(z \cdot y)l(x) + l(x)l(y \cdot z) - l(x)l(z \cdot y) - l(z \cdot x)l(y) + l(x \cdot z)l(y) + l(y)l(z \cdot x) - l(y)l(x \cdot z)$
\par
$= - [l(x),l(y)]l(z) + [l(y),l(x)]l(z) + l(z)[l(x),l(y)] - l(z)[l(y),l(x)] - [l(y),l(z)]l(x) + [l(z),l(y)]l(x) + l(x)[l(y),l(z)] - l(x)[l(z),l(y)] - [l(z),l(x)]l(y) + [l(x),l(z)]l(y) + l(y)[l(z),l(x)] - l(y)[l(x),l(z)]$ (again by (\ref{eq2.1}))
\par
$= - [[l(x),l(y)],l(z)] + [[l(y),l(x)],l(z)] + [l(x),[l(y),l(z)]] + [[l(z),l(y)],l(x)] + [l(y),[l(z),l(x)]]+ [[l(x),l(z)],l(y)] $
\par
$= - [[l(x),l(y)],l(z)] - [[l(y),l(z)],l(x)] - [[l(z),l(x)],l(y)] + [[l(y),l(x)],l(z)] + [[l(x),l(z)],l(y)] + [[l(z),l(y)],l(x)]$
\par
$= 0$\\
so (R2) holds for $(\rho, \theta, D)$.
\par
$\bullet$ For (R3) we have
\par
$\theta ([x,y],z) - \theta (x,z) \rho (y) + \theta (y,z) \rho (x)$
\par
$= -r(z)r(x \cdot y) + r(z)r(y \cdot x) + r(z)r(x)l(y) - r(z)r(x)r(y) - r(z)r(y)l(x) + r(z)r(y)r(x)$
\par
$= -r(z)\{ r(y)r(x) + l(x)r(y) \} + r(z)\{ r(x)r(y) + l(y)r(x) \} + r(z)r(x)l(y)-r(z)r(x)r(y) - r(z)r(y)l(x) + r(z)r(y)r(x)$ (by (\ref{eq2.2}))
\par
$ = -r(z)l(x)r(y) + r(z)r(x)l(y) + r(z)l(y)r(x)- r(z)r(y)l(x)$
\par
$= -r(z)l(x)r(y) - r(z)r(x)r(y) + r(z)l(y)r(x)+ r(z)r(y)r(x)$ (by (\ref{eq2.4}))
\par
$= -\{ r(z)l(x) + r(z)r(x) \}_{(2.4)} r(y) +
\{ r(z)l(y) + r(z)r(y) \}_{(2.4)} r(x)$
\par
$= 0$ (by (\ref{eq2.4}))\\
so (R3) holds for $(\rho, \theta, D)$.
\par
$\bullet$ For (R4) we have
\par
$[D(x,y), \rho (z)] = D(x,y)\rho (z) - \rho (z) D(x,y)$
\par
$= -l(x \cdot y)l(z) + l(x \cdot y)r(z) + l(z)l(x \cdot y) - r(z)l(x \cdot y)$
\par
$= - \{ l(x \cdot y)l(z) - l(z)l(x \cdot y)\}_{(2.1)} + l(x \cdot y)r(z)-r(z)l(x \cdot y)$
\par
$= -l((x \cdot y) \cdot z) + l(x \cdot y)r(z)-r(z)l(x \cdot y)$ (by (\ref{eq2.1}))
\par
$= -l((x \cdot y) \cdot z) + l(x \cdot y)r(z) +  r(z)r(x \cdot y)$ (by (\ref{eq2.4}))
\par
$= - l((x \cdot y) \cdot z) + r((x \cdot y) \cdot z)$ (by (\ref{eq2.2}))
\par
$= \rho (-(x \cdot y) \cdot z)$
\par
$= \rho ([\![x, y,z]\!])$\\
so (R4) holds for $(\rho, \theta, D)$.
\par
$\bullet$ For (R5) we have
\par
$\rho (y) \theta (x,z) - \rho (z) \theta (x,y) = - l(y)r(z)r(x) + r(y)r(z)r(x) + l(z)r(y)r(x) - r(z)r(y)r(x)$
\par
$= \{ -r(z)r(y) - l(y)r(z) \}_{(2.2)}r(x) + \{ r(y)r(z) + l(z)r(y) \}_{(2.2)} r(x)$
\par
$= -r(y \cdot z)r(x) + r(z \cdot y)r(x)$ (by (\ref{eq2.2}))
\par
$= -r([y,z])r(x)$
\par
$= \theta (x, [y,z])$\\
so we get (R5) for $(\rho, \theta, D)$.
\par
$\bullet$ For (R6) we have
\par
$\theta ([\![u, v,x]\!], y) + \theta (x, [\![u, v,y]\!]) = r(y)r((u \cdot v) \cdot x) + r((u \cdot v) \cdot y) r(x)$
\par
$= r(y) \{ r(x) r(u \cdot v) + l(u \cdot v) r(x)\} + \{ r(y) r(u \cdot v) + l(u \cdot v) r(y) \} r(x)$ (by (\ref{eq2.2}))
\par
$= r(y)r(x)r(u \cdot v) + l(u \cdot v) r(y) r(x) + \{ r(y)l(u \cdot v) + r(y)r(u \cdot v) \}_{(\ref{eq2.4})} r(x)$
\par
$= r(y)r(x)r(u \cdot v) + l(u \cdot v) r(y) r(x)$ (by (\ref{eq2.4}))
\par
$= - r(y)r(x)l(u \cdot v) + l(u \cdot v) r(y) r(x)$ (by (\ref{eq2.4}))
\par
$= [D(u,v), \theta (x,y)]$\\
and so (R6) holds for $(\rho, \theta, D)$.
\par
$\bullet$ Finally for (R7) we have
\par
$\theta (y,z) \theta (u,x) - \theta (x,z) \theta (u,y) + D(x,y) \theta (u,z)$
\par
$= r(z)r(y)r(x)r(u) - r(z)r(x)r(y)r(u) + l(x \cdot y)r(z)r(u)$
\par
$= r(z)(r(y)r(x) - r(x)r(y))r(u) + l(x \cdot y)r(z)r(u)$
\par
$= r(z)(r(y)r(x) + r(x)l(y))r(u) + l(x \cdot y)r(z)r(u)$ (by (\ref{eq2.4}))
\par
$= r(z)(r(x \cdot y) - l(x)r(y) + r(x)l(y))r(u)+ l(x \cdot y)r(z)r(u)$ (by (\ref{eq2.2}))
\par
$= r(z)r(x \cdot y)r(u) + r(z)(r(x)l(y) - l(x)r(y))r(u) + l(x \cdot y)r(z)r(u)$
\par
$= r(z)r(x \cdot y)r(u) + r(z)(-r(x)r(y) - l(x)r(y))r(u) + l(x \cdot y)r(z)r(u)$ (by (\ref{eq2.4}))
\par
$= r(z)r(x \cdot y)r(u) + \{ -r(z)r(x) - r(z)l(x) \}_{(\ref{eq2.4})}r(y)r(u) + l(x \cdot y)r(z)r(u)$
\par
$= r(z)r(x \cdot y)r(u)+l(x \cdot y)r(z)r(u)$ (by (\ref{eq2.4}))
\par
$= r((x \cdot y) \cdot z)r(u)$ (by (\ref{eq2.2}))
\par
$= \theta (u, [\![x,y,z]\!])$\\
so (R7) holds for $(\rho, \theta, D)$. This completes the proof.
\end{proof}
The notion of a homomorphism of representations of Leibniz algebras was introduced in \cite{LP}.\\

\noindent
{\bf Definition 2.7.} Two representations
$(V_1,l_1,r_1)$ and $(V_2,l_2,r_2)$ of a left Leibniz algebra $(\mathfrak{g}, \cdot)$ are said to be equivalent if there is an isomorphism $\psi : V_1 \rightarrow V_2$ such that $\psi l_1 (x) = l_2 (x) \psi$ and $ \psi r_1 (x) = r_2 (x) \psi$ for all $x \in \mathfrak{g}$.\\
\par
The notion of a homomorphism of representations of a $LY$ algebra was recently introduced in \cite{SZ}. From this notion one easily defines equivalent representations of a $LY$ algebra as follows.\\

\noindent
{\bf Definition 2.8.} Two representations
$(V_1, {\rho}_1 , {\theta}_1)$ and $(V_2, {\rho}_2 , {\theta}_2)$ of a given $LY$ algebra $\mathfrak{g}$ are said to be equivalent if there is an isomorphism $\psi : V_1 \rightarrow V_2$ such that $\psi {\rho}_1 (x) = {\rho}_2 (x) \psi$ and $\psi {\theta}_1 (x,y) = {\theta}_2 (x,y) \psi$ for all $x,y \in \mathfrak{g}$.\\
\par
The following result shows that equivalent representations of a given left Leibniz algebra give rise to equivalent representations of the associated $LY$ algebra.\\

\noindent
{\bf Theorem 2.9.} 
{\it Let $(\mathfrak{g}, \cdot)$ be a left Leibniz algebra and $(\mathfrak{g}, [,], [\![,,]\!])$ its associated $LY$ algebra. If
$(V_1,l_1,r_1)$ and $(V_2,l_2,r_2)$ are two equivalent representations of $(\mathfrak{g}, \cdot)$, then the associated representations
$(V_1, {\rho}_1 , {\theta}_1)$ and $(V_2, {\rho}_2 , {\theta}_2)$ of $(\mathfrak{g}, [,], [\![,,]\!])$ are also equivalent.}
\begin{proof}
We have, $\forall x,y \in \mathfrak{g}$,
$\forall v \in V_1$,
\par
$\psi ({\rho}_1 (x)(v)) = \psi (l_1 (x)(v)) - \psi (r_1 (x)(v))$
\par
$= l_2 (x)\psi (v) - r_2 (x)\psi (v)$
\par
$= {\rho}_2 (x) (\psi (v))$\\
and, from the other hand,
\par
${\theta}_2 (x,y) (\psi (v)) = - r_2 (y)(\psi (r_1 (x)(v))) = - \psi (r_1 (y)(r_1 (x)(v)))$
\par
$= \psi (- r_1 (y)(r_1 (x)(v))) = \psi ({\theta}_1 (x,y)(v))$\\
so that $\psi {\rho}_1 (x) = {\rho}_2 (x) \psi$ and $\psi {\theta}_1 (x,y) = {\theta}_2 (x,y) \psi$. Thus $(V_1, {\rho}_1 , {\theta}_1)$ and $(V_2, {\rho}_2 , {\theta}_2)$ are equivalent.
\end{proof}

\section{Examples}

In this section we consider some well-known representations of left Leibniz algebras, as well as their corresponding representations of associated Lie-Yamaguti algebras, after Theorem 2.6. Below $(\mathfrak{g}, \cdot)$ denotes a left Leibniz algebra and
$(\mathfrak{g}, [,], [\![,,]\!])$ its associated $LY$ algebra.\\

\noindent
{\bf Example 3.1.} The adjoint representation $(\mathfrak{g}, L, R)$ of $(\mathfrak{g}, \cdot)$ induces, by Theorem 2.6, a representation $(\mathfrak{g}, \rho, \theta)$ of  $(\mathfrak{g}, [,], [\![,,]\!])$, where
$\rho (x) := L(x) - R(x) = L_x - R_x$ and
$\theta (x,y) := -R(y)R(x) = -R_y R_x$,
$\forall x,y \in \mathfrak{g}$.
\par
Observe that $(\mathfrak{g}, \rho, \theta)$ is precisely the adjoint representation of
$(\mathfrak{g}, [,], [\![,,]\!])$. Thus the adjoint representation of $(\mathfrak{g}, \cdot)$ induces the one of $(\mathfrak{g}, [,], [\![,,]\!])$.\\

\noindent
{\bf Example 3.2.}  If $(V,l,r)$ is a left representation of $(\mathfrak{g}, \cdot)$, then $(V^*, l^*, -l^*-r^*)$ is a left representation of $(\mathfrak{g}, \cdot)$ \cite{TS} (see also Lemma 2.3 above). Therefore, Theorem 2.6 implies that $(V^* , \rho, \theta)$ is a representation of
$(\mathfrak{g}, [,], [\![,,]\!])$, where $\rho (x) := 2l^* (x) + r^* (x)$ and $\theta (x,y) := -(-l^* (y)-r^* (y))(- l^* (x)-r^* (x))$.
\par
Let $r_0 (z) := -l(z)-r(z)$, $\forall z \in \mathfrak{g}$. Then $(V^*, l^*, -l^*-r^*) = (V^*, l^*, r_0 ^*)$ and $\theta (x,y)=-r_0 ^* (y) r_0 ^* (x)$, $\rho (x) := l^* (x) - r_0 ^* (x)$.\\

\noindent
{\bf Example 3.3.} Suppose that
$(\mathfrak{g}, \cdot)$ is finite-dimensional and that $V$ is an irreducible finite-dimensional module for the Leibniz algebra
$(\mathfrak{g}, \cdot)$ (this is equivalent to say that the corresponding representation $l,r : \mathfrak{g} \rightarrow End(V)$ is irreducible). Then one knows (\cite{B1, B2, LP}) that either $r(x)=0$ or $r(x)=-l(x)$ for all $x \in \mathfrak{g}$. In \cite{LP} the representations $(V,l,r)$ of
 $(\mathfrak{g}, \cdot)$ for which $r(x)=0$,
 $\forall x \in \mathfrak{g}$, were said to be antisymmetric while the ones for which $r(x)=-l(x)$, $\forall x \in \mathfrak{g}$, were said to be symmetric.
 \par
 In case when $(V,l,r)$ is an antisymmetric representation of $(\mathfrak{g}, \cdot)$, Theorem 2.6 implies that the associated representation $(V, \rho, \theta)$ of
 $(\mathfrak{g}, [,], [\![,,]\!])$ is given by $\rho (x) = l(x)$, $\theta (x,y) = 0$, $\forall x,y \in \mathfrak{g}$.
 \par
 In case when $(V,l,r)$ is a symmetric representation of $(\mathfrak{g}, \cdot)$ then, by Theorem 2.6, the associated representation $(V, \rho, \theta)$ of
 $(\mathfrak{g}, [,], [\![,,]\!])$ is given by $\rho (x) = 2l(x)$, $\theta (x,y)=-l(y)l(x)$, $\forall x,y \in \mathfrak{g}$.

\end{document}